\pdfoutput=1
\documentclass[11pt,a4paper]{article}

\usepackage[margin=1in]{geometry}
\usepackage[utf8]{inputenc}
\usepackage{url}
\usepackage{booktabs}
\usepackage{amsfonts}
\usepackage{graphicx}
\usepackage{amsmath,amsfonts,amssymb,amsthm}
\usepackage{subcaption}
\usepackage{enumitem}
\usepackage{setspace}
\usepackage{lmodern}
\usepackage{caption}
\usepackage{xspace}
\usepackage{hyperref} 
\theoremstyle{plain}
\newtheorem{theorem}{Theorem}[section]
\newtheorem{lemma}[theorem]{Lemma}

\theoremstyle{definition}
\newtheorem{definition}[theorem]{Definition}

\theoremstyle{remark}
\newtheorem{remark}[theorem]{Remark}
\usepackage{titlesec}
\titlespacing*{\section}{0pt}{6pt}{4pt}
\titlespacing*{\subsection}{0pt}{4pt}{3pt}

\newcommand{\Poincare}{Poincar\'e\xspace}

\def\keywordname{{\bfseries \emph{Keywords}}}%
\def\keywords#1{\par\addvspace\medskipamount{\rightskip=0pt plus1cm
\def\and{\ifhmode\unskip\nobreak\fi\ $\cdot$
}\noindent\keywordname\enspace\ignorespaces#1\par}}

\newtheoremstyle{plain}
  {3pt} 
  {3pt}
  {\itshape}
  {}
  {\bfseries}
  {.}
  {0.5em}
  {}
\theoremstyle{plain}

\date{}
\title{Synchronization Phenomenon in Three-Time-Scale Systems}
\author{
  \href{https://orcid.org/0000-0003-4207-3032}{Navojit Dhali Pallab} \\
  Mathematical Institute, Tohoku University, Sendai,980-8578, Japan\\
  \texttt{pallab.navojit.dhali.t4@dc.tohoku.ac.jp}
}

\begin{document}
\normalfont
\maketitle

\begin{abstract}
This paper investigates synchronization phenomena in networks of coupled oscillators governed by three-time-scale dynamical systems exhibiting canard dynamics. A mathematical framework has been developed to analyze the synchronization of fast variables across heterogeneous systems, deriving a sufficient condition for the synchronization error to fall below a specified threshold within the minimum linger time. This condition accounts for coupling strength, heterogeneity, and time-scale separation, ensuring stable oscillatory behavior in the network. The result, supported by rigorous mathematical analysis, advances the understanding of synchronization in complex multi-time-scale systems. 
\end{abstract}

\def\keywordname{\textbf{Keywords:}}
\def\keywords#1{\par\addvspace{0.5\baselineskip}\noindent\keywordname\ \textit{#1}\par\vspace{0.5\baselineskip}}

\section{Introduction}
Synchronization is a fundamental phenomenon in the network of oscillators, particularly when the units within the network exhibit heterogeneity. In the real world, systems such as neuronal networks, cardiac tissues, or the pancreatic $\beta-$cells network, these differences arise naturally due to variations in cellular properties. For example, neurons may have varying firing thresholds, cardiac cells may differ in excitability, and $\beta$-cells may exhibit variations in glucose sensitivity or ion channel properties \cite{gutierrez2017heterogeneity, schuit1996factors}. Unlike homogeneous systems, where identical units can perfectly align with each other, heterogeneity complicates synchronization, often preventing exact synchrony and making the study of near-synchronization a rich and challenging topic.

In our study, synchronization does not necessarily imply identical behavior across all units. Instead, we define synchronization as a state where the differences between key variables—specifi-cally, the fast variables $v_i$ (membrane potentials in $\beta$-cells) remain small or nearly identical within a tolerance $\epsilon$. This relaxed definition, measured via the synchronization error $V_v = \frac{1}{N} \sum_{i=1}^N ( \bar{v}-v_i)^2$, is practical for biological systems, where exact synchrony is often neither feasible nor necessary. The physiological importance of such synchronization is profound. Coordinated neuronal firing underpins cognitive functions like memory and perception in the brain. In the heart, synchronized contractions ensure efficient blood pumping. In pancreatic $\beta$-cell networks, synchronized bursting facilitates the coordinated release of insulin to regulate blood glucose levels \cite{sherman1991model}. Understanding the conditions under which synchronization occurs in heterogeneous networks is thus critical for both theoretical insights and practical applications, such as diagnosing and treating disorders like epilepsy, arrhythmias, or diabetes-related dysfunctions.

This chapter focuses on synchronization in networks of coupled oscillators with three distinct time scales—fast, intermediate, and slow—exhibiting canard dynamics, a feature common in biological systems \cite{wechselberger2020geometric, kuehn2015multiple}. While synchronization in coupled oscillators has been extensively studied, most analyses focus on two-time-scale systems \cite{strogatz2000kuramoto, pikovsky2001synchronization,chen2018synchronization}, leaving a gap in understanding three-time-scale dynamics or multi-time-scale dynamics. Pikovsky et al. \cite{pikovsky2001synchronization} developed a comprehensive framework for heterogeneous systems. Building on this foundation, we derive a sufficient condition for synchronization, ensuring that the synchronization error $V_v$ falls below the tolerance $\epsilon$ within a specified time window $t_{linger}^{\min}$, the minimum linger time across the network accounting for variability due to heterogeneity. The general condition is
\[
k > \max \left( \frac{2M}{\sqrt{\epsilon}}, \frac{1}{\delta t_{linger}^{\min}} \ln \left( \frac{2 W_0}{\sqrt{\epsilon}} \right) \right),
\]
where $k$ is the coupling strength, $M$ bounds the heterogeneity, $W_0$ bounds the initial synchronization error, and $\delta \in (0, 1)$ is the slowest time scale parameter. In near-homogeneous cases with specific parameter choices (e.g., $M \sim 0$, $\delta \approx 1$, $W_0 = \frac{1}{2}$), this condition simplifies to forms like $k > -\frac{\ln \epsilon}{2 t_{linger}^{\min}}$.

The paper has the following two sections: (i) presents the model of a three-time-scale network and statement of the main theorem (Section \ref{sec:mainresultsTheory}), and (ii) analyzes the synchronization error through a series of mathematical derivations with proof of the theorem (Section \ref{sec:analysisANDproof}).

\section{Main Theorem}
\label{sec:mainresultsTheory}
In this section, we state our main theorem about synchronization in a network of $N$ coupled three-time-scale dynamical systems, each evolving in $\mathbb{R}^5$. These systems model phenomena with distinct fast, intermediate, and slow dynamics, such as slow bursters in biological networks (e.g., $\beta-$cells network), where time scale separation is critical in the emergence of complex oscillatory behavior. The following system of differential equations defines the network
\begin{equation}
\begin{aligned}
\frac{dv_i}{dt} &= h_{1i}(v_i, u_i, x_i, y_i, z_i, \varepsilon, \delta; \mu_i) + k \sum_{j=1}^N a_{ij} (v_j - v_i), \\
\frac{du_i}{dt} &= h_{2i}(v_i, u_i, x_i, y_i, z_i, \varepsilon, \delta; \mu_i), \\
\frac{dx_i}{dt} &= \varepsilon f_i(v_i, u_i, x_i, y_i, z_i, \varepsilon, \delta; \mu_i), \\
\frac{dy_i}{dt} &= \varepsilon \delta g_{1i}(v_i, u_i, x_i, y_i, z_i, \varepsilon, \delta; \mu_i), \\
\frac{dz_i}{dt} &= \varepsilon \delta g_{2i}(v_i, u_i, x_i, y_i, z_i, \varepsilon, \delta; \mu_i),
\end{aligned}\label{sys:general3timescale5DNetwork}
\end{equation}
for $i = 1, \ldots, N$. Here, $\mathbf{x}_i = (v_i, u_i, x_i, y_i, z_i) \in \mathbb{R}^5$ represents the state of the $i$-th system, where $v_i, u_i \in \mathbb{R}$ are fast variables, $x_i \in \mathbb{R}$ is an intermediate variable, and $y_i, z_i \in \mathbb{R}$ are slow variables. The small positive constants $\varepsilon, \delta $ define the time scale separation ($0<\varepsilon, \delta \ll 1$), with $\varepsilon \delta$ governing the slowest dynamics, while $\mu_i \in \mathbb{R}^m$ is a control parameter (e.g., representing glucose concentration, KATP conductance, glucose sensitivity or, up-scaling in $\beta$-cell models). The functions $h_{1i}, h_{2i}, f_i, g_{1i}, g_{2i}: \mathbb{R}^5 \times (0, 1) \times (0, 1) \times \mathbb{R}^m \to \mathbb{R}$ are assumed to be smooth ($C^\infty$). The coupling strength $k \geq 0$ governs the interaction between systems, with coupling coefficients $a_{ij} = \frac{1}{N}$, corresponding to all-to-all coupling satisfying $\sum_{j=1}^N a_{ij} = 1$.

When $N = 1$, the system \eqref{sys:general3timescale5DNetwork} reduces to the uncoupled system. Moreover, for each $i$, when the coupling strength $k = 0$, the system \eqref{sys:general3timescale5DNetwork} inherits the properties of the uncoupled system, including the existence of critical manifolds, canard points, and slow bursting behavior. We only provide the essential background on the system's structure and dynamics to maintain focus on synchronization.

\subsection{Key definition and Theorem}
\begin{definition}[Linger Time $t_{linger}^i$ and Synchronization Time Window]
\label{def:tlinger} 
Consider the network of $N$ three-time-scale dynamical systems defined by Equation \eqref{sys:general3timescale5DNetwork}, with fast variables $(v_i, u_i)$, intermediate variable $x_i$, slow variables $(y_i, z_i)$, and time scale parameters $0 < \varepsilon, \delta \ll 1$. Due to heterogeneity, the intrinsic dynamics $(h_{1i}, h_{2i}, f_i, g_{1i}, g_{2i})$ vary slightly across systems. For each system $i$ (with coupling $k = 0$), the fast critical manifold $S_i$ is 
\begin{align*}\
    S_i = \{ (v_i, u_i, x_i, y_i, z_i) \in \mathbb{R}^5 \mid h_{1i}(v_i, u_i, x_i, y_i, z_i, 0, \delta; \mu_i)= 0,\ h_{2i}(v_i, u_i, x_i, y_i, z_i, 0, \delta; \mu_i) = 0 \},
\end{align*}
where $S_i$ admits a local parameterization $v_i = \phi_v^i(x_i, y_i, z_i; \delta, \mu_i)$, $u_i = \phi_u^i(x_i, y_i, z_i; \delta, \mu_i)$, and consists of an attracting branch $S_{a,i}$, a repelling branch $S_{r,i}$, and a periodic branch $S_{p,i}$. The fold curve of $S_i$, where trajectories jump from $S_{a,i}$ to $S_{p,i}$, is given by the conditions $h_{1i} = h_{2i} = 0$ and $\det(D_{(v_i,u_i)} \mathbf{h}_i) = 0$, with a jump point denoted $(v_i^f, u_i^f, x_i^f, y_i^f, z_i^f)$.

The slow critical manifold $M_i$, a submanifold of $S_i$, is
\begin{align*}
    M_i = \{ (v_i, u_i, x_i, y_i, z_i) \in S_i \mid f_i(v_i, u_i, x_i, y_i, z_i, 0, 0; \mu_i) = 0 \}.
\end{align*}
Assuming $\partial f_i / \partial x_i \neq 0$ on $S_i$, the implicit function theorem ensures that $M_i$ can be locally parameterized by the slow variables $y_i$ and $z_i$, yielding functions $\psi_v^i(y_i, z_i; \mu_i)$, $\psi_u^i(y_i, z_i; \mu_i)$, and $\psi_x^i(y_i, z_i; \mu_i)$ such that
\begin{align*}
    v_i = \psi_v^i(y_i, z_i; \mu_i), \ \  u_i = \psi_x^i(y_i, z_i; \mu_i), \ \ x = \psi_x^i(y_i, z_i; \mu_i),
\end{align*}
for $(v_i, u_i, x_i, y_i, z_i) \in M_i$. A canard point $(v_i^c, u_i^c, x_i^c, y_i^c, z_i^c) \in S_{a,i} \cap M_i$ marks the transition where the trajectory of system $i$ lingers on $S_{a,i} \cap M_i$ before jumping to $S_{p,i}$.

We define the entry \Poincare section near the canard point as
\begin{align*}
\Sigma_{in,i}^c = \{ (v_i, u_i, x_i, y_i, z_i) \mid v_i = \phi_v^i(x_i^c - \delta_x, y_i^c, z_i^c; \delta, \mu_i),\ u_i = \phi_u^i(x_i^c - \delta_x, y_i^c, z_i^c; \delta, \mu),\\ x_i = x_i^c - \delta_x,\ |y_i - y_i^c| < \delta_y,\ |z_i - z_i^c| < \delta_z \},
\end{align*}
and the pre-jump \Poincare section as

\begin{align*}
\Sigma_{pre-jump,i} = \{ (v_i, u_i, x_i, y_i, z_i) \mid v_i = \phi_v^i(x_i^f - \delta_x', y_i^f, z_i^f; \delta, \mu_i),\ u_i = \phi_u^i(x_i^f - \delta_x', y_i^f, z_i^f; \delta, \mu_i), \\x_i = x_i^f - \delta_x', \ |y_i - y_i^f| < \delta_y',\ |z_i - z_i^f| < \delta_z' \},
\end{align*}
where $\delta_x, \delta_x', \delta_y, \delta_y', \delta_z, \delta_z' > 0$ are small constants. The \textit{linger time} $t_{\text{linger}}^i$ for system $i$ is the duration its trajectory spends on $S_{a,i} \cap M_i$, from $\Sigma_{in,i}^c$ to $\Sigma_{\text{pre-jump},i}$, determined by the slowest dynamics in $y$ (one of the slowest variables)
\begin{equation}
\begin{aligned}
t_{\text{linger}}^i = \int_{y_i^c - \delta_y}^{y_i^f} \frac{dy}{\varepsilon \delta \tilde{g}_1^i(y, z_i^c;  \mu_i)},
\end{aligned}
\end{equation}
where
\begin{align*}
    \tilde{g}_1^i(y, z_i^c;\mu) = g_{1i}(\psi_v^i(y, z_i^c;  \mu_i), \psi_u^i(y, z_i^c;\mu_i), \psi_x^i(y, z_i^c; \mu_i), y, z_i^c, 0, 0; \mu_i).
\end{align*}

Because of the heterogeneity, $t_{linger}^i$ varies slightly across systems. For synchronization purposes, define the \textit{synchronization time window} as
\begin{align*}
    t_{linger}^{\min} = \min_{i=1,\ldots,N} \{ t_{linger}^i \},
\end{align*}
ensuring all systems are still in the lingering phase when synchronization is achieved.
\end{definition}

\begin{theorem}
\label{thm:synchronization}
Consider a network of $N$ three-time-scale dynamical systems in $\mathbb{R}^5$, defined by Equation \ref{sys:general3timescale5DNetwork}, with coupling matrix $a_{ij} = \frac{1}{N}$. Let the fast ($\varepsilon=0$) critical manifold $S_i$ for each system $i$ have attracting branch $S_{a,i}$, repelling branch $S_{r,i}$, and periodic branch $S_{p,i}$, and let $M_i$ be the slow ($\varepsilon=\delta=0$) critical manifold. Assume 
\begin{enumerate}
    \item[(i)]  The intrinsic dynamics satisfy $|h_{1i}(v_i, u_i, x_i, y_i, z_i, \varepsilon, \delta; \mu_i)| \leq M$ for all $i$, where $M \geq 0$ is a constant bounding the heterogeneity (Remark~\ref{remarkonM} on $M$).
    \item[(ii)] For each system $i$, there exists a canard point $(v_i^c, u_i^c, x_i^c, y_i^c, z_i^c) \in S_{a,i} \cap M_i$, and a jump point on the fold curve of $S_i$, such that the trajectory lingers on $S_{a,i} \cap M_i$ before transitioning to $S_{p,i}$.
    \item[(iii)] The linger time $t_{\text{linger}}^i$, as defined in Definition \ref{def:tlinger}, varies across systems due to heterogeneity, and the synchronization time window is $t_{linger}^{\min}=\min_i\{ t_{linger}^i\}$.
    \item[(iv)] $v_i(t)\in S_{i,a}$ (stable branch of $S_i$), for all $i=1,\cdots,N$, for $t\in[0,T]$ with  $t_{linger}^{\min} \leq \frac{T}{\delta}$.
\end{enumerate}
Let the initial synchronization error satisfy $W(0) = \sqrt{V_v(0)}\leq W_0$ for some positive constant $W_0$, where
\begin{align*}
    V_v(t) = \frac{1}{N} \sum_{i=1}^N (v_i(t)-\bar{v}(t))^2, \ \ \ \bar{v}(t)=\frac{1}{N}\sum_{j=1}^N v_j(t).
\end{align*}
A \emph{sufficient condition} for synchronization, with $V_v(t_{linger}^{\min})<\epsilon$, is that the coupling strength $k$ satisfies
\begin{equation}
\begin{aligned}
k > \max \left( \frac{2M}{\sqrt{\epsilon}}, \frac{1}{\delta t_{linger}^{\min}} \ln \left(\frac{2 W_0}{\sqrt{\epsilon}} \right) \right),
\end{aligned}
\end{equation}
where $\epsilon > 0$ is the desired synchronization error threshold, and $\delta \in (0, 1)$ is the slowest time scale separation from the intermediate time scale $\varepsilon$.
\end{theorem}

\begin{remark}[Sufficiency and Practical Significance]
\label{rem:sufficiency_significance}
The condition in Theorem \ref{thm:synchronization} is sufficient but not necessary, meaning synchronization may occur for smaller $k$ if the intrinsic dynamics (attraction to $S_i$) or lower heterogeneity assists the process. This sufficient condition is valuable as it provides a practical threshold for ensuring synchronization in complex three-time-scale systems, particularly in applications like $\beta$-cell networks, where the predicted coupling strengths match physiological observations. The adjustment $\delta t_{linger}$ accounts for the slowest time scale, offering a different perspective on synchronization in such systems.
\end{remark}

\section{Analysis and Proof}\label{sec:analysisANDproof}
To prove Theorem \ref{thm:synchronization}, we have to first analyze the dynamics of the variance $V_v$ of the fast variable $v_i$, defined as
\begin{align}\label{eq:sycError}
    V_v(t) = \frac{1}{N} \sum_{i=1}^N (v_i(t) - \bar{v}(t))^2, \quad \bar{v}(t) = \frac{1}{N} \sum_{j=1}^N v_j(t).
\end{align}

This variance-based metric is a standard measure of synchronization in coupled dynamical systems, commonly used in oscillator networks and biological systems \cite{arenas2008synchronization,pikovsky2001synchronization}. Our goal is to derive the time evolution of $V_v$ and identify conditions under which $V_v(t_{linger}) < \epsilon$, indicating near-synchronous behavior of the $v_i$-variables across the network.

Consider the network of $N$ three-time-scale dynamical systems defined by Equation (\ref{sys:general3timescale5DNetwork}). The dynamics of each system $i$ are
\begin{align*}
    \dot{\mathbf{x}}_i = \mathbf{f}_i(\mathbf{x}_i) + k \sum_{j=1}^N a_{ij} (v_j - v_i) \mathbf{e}_1,
\end{align*}
where $\mathbf{x}_i = (v_i, u_i, x_i, y_i, z_i)^T \in\mathbb{R}^5$. $\mathbf{f}_i(\mathbf{x}_i) = (h_{1i}, h_{2i}, \varepsilon f_i, \varepsilon \delta g_{1i}, \varepsilon \delta g_{2i})^T$, and $\mathbf{e}_1\ = (1, 0, 0, 0, 0)^T$. Through this section, $\dot{X}$ represents the derivative of $X$ with respect to $t$, where $t$ is the fast time scale. The diffusive coupling term, inspired by the Kuramoto model \cite{acebron2005kuramoto, kuramoto1984chemical}, drives the $v_i$-variables toward their mean. With all-to-all coupling $a_{ij} = \frac{1}{N}$, the coupling term simplifies to
\begin{align*}
    k \sum_{j=1}^N a_{ij} (v_j - v_i) = k \left( \bar{v} - v_i \right).
\end{align*}
Thus, the dynamics of the $v_i$-component are
\begin{equation}
\begin{aligned}
\dot{v}_i = h_{1i}(v_i, u_i, x_i, y_i, z_i, \varepsilon, \delta; \mu) + k (\bar{v} - v_i).
\end{aligned}
\label{eq:v_i_dynamics}
\end{equation}

\begin{lemma}
For $N-$coupled systems (\ref{sys:general3timescale5DNetwork}), the variance metric (\ref{eq:sycError}) satisfied 
    \begin{equation*}
    \begin{aligned}
        \dot{V}_v = -2k V_v + \frac{2}{N} \sum_{i=1}^N (v_i - \bar{v}) (h_{1i} - \bar{h}_1).
    \end{aligned}
    \end{equation*}
\end{lemma}

\begin{proof}
    The derivative of $V_v$ with respect to $t$ in equation (\ref{eq:sycError}), following a standard approach in synchronization studies, gives
    \begin{align*}
        \dot{V}_v = \frac{d}{dt} \left( \frac{1}{N} \sum_{i=1}^N (v_i - \bar{v})^2 \right) = \frac{2}{N} \sum_{i=1}^N (v_i - \bar{v}) \frac{d}{dt} (v_i - \bar{v}) = \frac{2}{N} \sum_{i=1}^N (v_i - \bar{v}) (\dot{v}_i - \dot{\bar{v}}).
    \end{align*}
    Using (\ref{eq:v_i_dynamics}), for the dynamics of the mean $\bar{v}$, we obtain
    \begin{align*}
        \dot{\bar{v}} = \frac{1}{N} \sum_{i=1}^N \dot{v}_i = \frac{1}{N} \sum_{i=1}^N \left[ h_{1i}(v_i, u_i, x_i, y_i, z_i, \varepsilon, \delta; \mu_i) + k (\bar{v} - v_i) \right].
    \end{align*}
    The coupling term vanishes as
    \begin{align*}
        \frac{k}{N} \sum_{i=1}^N (\bar{v} - v_i) = \frac{k}{N} (N \bar{v} - N \bar{v}) = 0.
    \end{align*}
    Therefore,
    \begin{align*}
        \dot{\bar{v}} = \frac{1}{N} \sum_{i=1}^N h_{1i}(v_i, u_i, x_i, y_i, z_i, \varepsilon, \delta; \mu_i) = \bar{h}_1.
    \end{align*}
    Substituting into the expression for $\dot{V}_v$, we have
    \begin{align*}
        \dot{V}_v = \frac{2}{N} \sum_{i=1}^N (v_i - \bar{v}) \left[ h_{1i}(v_i, u_i, x_i, y_i, z_i, \varepsilon, \delta; \mu_i) + k (\bar{v} - v_i) - \bar{h}_1 \right].
    \end{align*}
    This can be split into two terms
    \begin{align*}
        \dot{V}_v = \frac{2}{N} \sum_{i=1}^N (v_i - \bar{v}) (h_{1i}(v_i, u_i, x_i, y_i, z_i, \varepsilon, \delta; \mu_i) - \bar{h}_1) - \frac{2k}{N} \sum_{i=1}^N (v_i - \bar{v})^2.
    \end{align*}
    Since $\frac{1}{N} \sum_{i=1}^N (v_i - \bar{v})^2 = V_v$, the coupling term becomes $-2k V_v$, leading to
    \begin{equation}
    \begin{aligned}
    \dot{V}_v = -2k V_v + \frac{2}{N}  \sum_{i=1}^N (v_i - \bar{v}) (h_{1i}(v_i, u_i, x_i, y_i, z_i, \varepsilon, \delta; \mu_i) - \bar{h}_1).
    \end{aligned}
    \label{eq:V_v_dynamics}
    \end{equation}
\end{proof}

The second term in (\ref{eq:V_v_dynamics}) reflects the effect of heterogeneity in the intrinsic dynamics, while the first term shows that coupling reduces the synchronization error at a rate proportional to $V_v$, a phenomenon well-documented in synchronization literature \cite{arenas2008synchronization}.

To proceed, we bound the heterogeneity term in Equation \eqref{eq:V_v_dynamics}. Assume the intrinsic dynamics are bounded, a common assumption in synchronization studies to control the effect of heterogeneity \cite{boccaletti2006complex}, with
\begin{align*}
    |h_{1i}(v_i, u_i, x_i, y_i, z_i, \varepsilon, \delta; \mu_i)| \leq M \quad \text{for all } i,
\end{align*}
where $M \geq 0$ is a constant. Since $|\bar{h}_1| \leq M$, the difference is bounded by
\begin{align}\label{eq:heteroBound}
    |h_{1i}(v_i, u_i, x_i, y_i, z_i, \varepsilon, \delta; \mu_i) - \bar{h}_1| \leq |h_{1i}(v_i, u_i, x_i, y_i, z_i, \varepsilon, \delta; \mu_i)| + |\bar{h}_1| \leq 2M.
\end{align}

\begin{lemma}\label{lem:v_inequa_dyna}
    The evolution of $V_v$ satisfies the nonlinear differential inequality
    \begin{equation*}
    \begin{aligned}
    \dot{V}_v \leq -2k V_v + 4M \sqrt{V_v}.
    \end{aligned}
    \end{equation*}
\end{lemma}
\begin{proof}
    Applying the Cauchy-Schwarz inequality to the intrinsic dynamic terms in (\ref{eq:V_v_dynamics}), we obtain
    \begin{align*}
        &\left| \sum_{i=1}^N (v_i - \bar{v}) (h_{1i}(v_i, u_i, x_i, y_i, z_i, \varepsilon, \delta; \mu_i) - \bar{h}_1) \right| \\& \leq \sqrt{ \sum_{i=1}^N (v_i - \bar{v})^2 } \sqrt{ \sum_{i=1}^N (h_{1i}(v_i, u_i, x_i, y_i, z_i, \varepsilon, \delta; \mu_i) - \bar{h}_1)^2 }.
    \end{align*}
    Since $\sum_{i=1}^N (v_i - \bar{v})^2 = N V_v$, and (\ref{eq:heteroBound}) implies $\sum_{i=1}^N (h_{1i}(v_i, u_i, x_i, y_i, z_i, \varepsilon, \delta; \mu_i) - \bar{h}_1)^2 \leq 4M^2 N$, we obtain
    \begin{align*}
        \left| \frac{2}{N} \sum_{i=1}^N (v_i - \bar{v}) (h_{1i}(v_i, u_i, x_i, y_i, z_i, \varepsilon, \delta; \mu_i) - \bar{h}_1) \right| \leq 4M \sqrt{V_v}.
    \end{align*}
    Thus,
    \begin{equation}
    \begin{aligned}
    \dot{V}_v \leq -2k V_v + 4M \sqrt{V_v}.
    \end{aligned}
    \label{eq:V_v_inequality}
    \end{equation}
\end{proof}
The inequality (\ref{eq:V_v_inequality}) captures the balance between the stabilizing effect of coupling ($-2k V_v$) and the destabilizing effect of heterogeneity ($4M \sqrt{V_v}$).

\begin{remark}\label{remarkonM}
 For the $i-$th system, the heterogeneity term $h_{1i}(v_i, u_i, x_i, y_i, z_i, \varepsilon, \delta; \mu_i)$ governs the intrinsic dynamics of the fast variable for a single system. When it's uncoupled ($k=0$), the system satisfies the following
\begin{itemize}
    \item[(i)] In the case $M$ close to $0$, the dynamics of the fast variable $v_i$ are on the slow manifold $S_i$, more specifically, near the stable branch of the critical manifold for small $\varepsilon$. Mathematically, this is the consequence of the fact that the dynamics of the system on the critical manifold $S_i$ is governed by 
    \begin{align*}
        \frac{dx_i}{d\tau}= g_{i}(v_i, u_i, x_i, y_i, z_i, 0, \delta; \mu_i),
    \end{align*}
    where $\tau=\varepsilon t$ represents the intermediate time scale. For small $M$, $\mid S_i -S_j \mid <\eta$, where $\eta$ is a small positive constant depends on $M$ and $S_i$ and $S_j$ is the corresponding critical manifold of the $i^{th}$ and $j^{th}$ systems, respectively. It implies that the individual (uncoupled) dynamics of the systems are in the quiescent phase, which happens only when the fast dynamics is near (by the Fenichel theory \cite{fenichel1979geometric}, $\varepsilon-$close) to the normally hyperbolic attracting branch of $S_i$.
    \item[(ii)] For the case $M=0$, implies $\mid S_i-S_j \mid=0$, the critical manifold and there fold curve are identical.
\end{itemize}
\end{remark}

To analyze the inequality \ref{eq:V_v_inequality}, we linearize it for small $V_v$, a technique often used in synchronization studies to simplify nonlinear dynamics. Define $W = \sqrt{V_v}$, so $V_v = W^2$ and $\dot{V}_v = 2W \dot{W}$. Substituting into \ref{eq:V_v_inequality}, we obtain
\begin{align*}
    2W \dot{W} \leq -2k W^2 + 4M W.
\end{align*}
For $W > 0$, the inequality becomes
\begin{equation}
\begin{aligned}
\dot{W} \leq -k W + 2M.
\end{aligned}
\label{eq:W_inequality}
\end{equation}
This linear differential inequality describes the evolution of $W$, a desynchronization measure, where the term $-k W$ promotes synchronization, and $2M$ reflects the residual error due to heterogeneity.
\begin{lemma}\label{lem:W_inequality}
Consider that all the $i=1,\cdots,N$ remains on the stable branch of the critical manifold $S_i$ for a finite time interval $[0,T]$. Then,
    \begin{equation*}
    \begin{aligned}
    W(t) \leq \left( W(0) - \frac{2M}{k} \right) e^{-kt} + \frac{2M}{k}.
    \end{aligned}
    \end{equation*}
for $t\in[0,T]$, with $W(0)<W_0$.
\end{lemma}

\begin{lemma}[Gronwall Lemma]\label{GronwallLemma}
Let $\eta(t)\in C^1([0,T];\mathbb{R}^+)$ satisfy the differential inequality
        \begin{equation*}
            \begin{aligned}
                \frac{d\eta(t)}{dt} \leq a \eta(t) + \psi (t), \hspace{10mm} \eta(0)=\eta_0,
            \end{aligned}
        \end{equation*}
        where $a\in \mathbb{R}$ and $\psi(t) \in L^1 ([0,T];\mathbb{R}^+)$. Then 
        \begin{equation*}
            \begin{aligned}
                \eta(t) \leq \exp{(at)}(\eta_0 + \int_0^t \exp{(-as)} \psi (s)ds),
            \end{aligned}
        \end{equation*}
        for all $t \in [0,T]$.
\end{lemma}

\begin{proof}[Proof of Lemma \ref{lem:W_inequality}]
    Since on $S_{i,a}$, the dynamics of the system \ref{sys:general3timescale5DNetwork} governs by the slow variables, $W(t) \in C^1([0,T])$. Now, using Gronwall lemma (Lemma \ref{GronwallLemma}), we can easily obtain that the solution to the inequality \eqref{eq:W_inequality} satisfies
    \begin{equation}
    \begin{aligned}
    W(t) \leq \left( W(0) - \frac{2M}{k} \right) e^{-kt} + \frac{2M}{k}.
    \end{aligned}
    \label{eq:W_bound}
    \end{equation}
\end{proof}
The solution comprises a transient term, $\left( W(0) - \frac{2M}{k} \right) e^{-kt}$, which decays exponentially at rate $k$, and a steady-state term, $\frac{2M}{k}$, representing the residual error due to heterogeneity as $t \to \infty$. For large $k$, the steady-state error $\frac{2M}{k}$ becomes small, confirming that strong coupling enhances synchronization, consistent with findings in oscillator networks \cite{pikovsky2001synchronization}.

\subsection{Synchronization Condition with Linger-Time}
\label{subsec:sync_condition}

Our objective is to ensure synchronization within the synchronization time window $t_{linger}^{\min} = \min_i \{ t_{linger}^i \}$, defined in Definition \ref{def:tlinger}, such that the synchronization error satisfies $V_v(t_{linger}^{\min}) < \epsilon$. This choice ensures all systems are still in the lingering phase, as $t_{linger}^{\min}$ is the earliest time any system may jump to $S_{p,i}$. Since $V_v = W^2$, the condition is equivalent to
\begin{align*}
    W(t_{\text{linger}}^{\min}) < \sqrt{\epsilon}.
\end{align*} 

\begin{proof}[Proof of Theorem \ref{thm:synchronization}]
Assume the conditions of Theorem \ref{thm:synchronization} hold. From Lemma \ref{lem:v_inequa_dyna}, we have the inequality
\begin{align*}
    \dot{V}_v \leq -2k V_v + 4M \sqrt{V_v}.
\end{align*}
From Lemma \ref{lem:W_inequality}, for $t\in [0,T]$, linearization with $W = \sqrt{V_v}$, yields
\begin{align*}
    W(t) \leq \left( W(0) - \frac{2M}{k} \right) e^{-kt} + \frac{2M}{k}.
\end{align*}

Now, we have to separate the steady-state and transient contributions.

Consider the residual error to be less than the tolerance
\begin{align*}
    \frac{2M}{k} < \sqrt{\epsilon},
\end{align*}
which implies 
\begin{align*}
    k > \frac{2M}{\sqrt{\epsilon}}.
\end{align*}
This ensures that, in the long term, the error due to heterogeneity is sufficiently small. Here, the steady-state term $\frac{2M}{k}$ is independent of the time window. 

The transient term must decay quickly enough to achieve synchronization in the desired time window. If $W(0)>\frac{2M}{k}$, we have to bound the transient contribution to $\frac{\sqrt{\epsilon}}{2}$, and obtain
\begin{align*}
    \left( W(0) - \frac{2M}{k} \right) e^{-k t} < \frac{\sqrt{\epsilon}}{2},
\end{align*}
for $t\in[0,T]$. Assumption (iv) in Theorem \ref{thm:synchronization}, ensures that $\delta t_{linger}^{\min}<T$. Therefore,
\begin{align*}
    \left( W(0) - \frac{2M}{k} \right) e^{-k \delta t_{linger}^{\min}} < \frac{\sqrt{\epsilon}}{2}.
\end{align*}
Solving for $k$, and approximating for large $k$, $\frac{2M}{k}\ll W(0)$, we get
\begin{align*}
    k > \frac{1}{\delta t_{linger}^{\min}} \ln \left( \frac{2 W_0}{\sqrt{\epsilon}} \right).
\end{align*}

Combining this with the residual error, the condition becomes
\begin{align*}
    k > \max \left( \frac{2M}{\sqrt{\epsilon}}, \frac{1}{\delta t_{linger}^{\min}} \ln \left( \frac{2 W_0}{\sqrt{\epsilon}} \right) \right).
\end{align*}
This condition ensures synchronization by the earliest linger time, accounting for variability due to heterogeneity.
\end{proof}

\subsection{Remark on Practical Implication}
    This condition balances the steady-state requirement (controlling the residual error due to heterogeneity) and the transient requirement (ensuring the rapid decay of the initial error within the slowest time-scale-adjusted window $\delta t_{linger}^{\min}$). The framework provides a practical tool for designing coupling strengths to achieve synchronization in heterogeneous oscillator networks, with direct applications to biological systems, such as the $\beta-$cell network \cite{pallab2025pancreatic}, where coordinated bursting is essential for physiological function.

\section*{Acknowledgments}
This work forms part of the ongoing PhD thesis ``Canard Dynamics and Synchronization in the Network of Three-Time-Scale Systems'' at Tohoku University. The author gratefully acknowledges Professor Hayato Chiba (AIMR, Tohoku University, Japan) for his invaluable guidance, continuous support, and patience throughout this research. This work was supported by JST University Fellowships for Science and Technology Innovation (Grant No. JPMJFS2102) and JST SPRING (Grant No. JPMJSP2114).


\end{document}